\input amstex
\magnification=\magstep1
\documentstyle{amsppt}
\pagewidth{5.3in}
\pageheight{7.2in}
\define\il{\int\limits_}
\define\hvi{\varphi}
\define\eps{\varepsilon}
\define\Hvi{\Phi}
\define\ovr{\overline}

\define\om{\omega}

\define\al{\alpha}
\define\gm{\gamma}
\define\th{\theta}

\define\Gm{\Gamma}

\define\Dl{\Delta}
\define\dl{\delta}
\define\opr{\operatorname}
\define\bd{\partial}

\define\sm{\setminus}
\define\sbs{\subset}
\define\sps{\supset}

\define\Sbs{\Subset}

\define\cl{\operatorname{cl}}
\define\supp{\operatorname{supp}}
\define\ess{\operatorname{ess}}
\define\dist{\operatorname{dist}}
\define\bl{\bold}
\document
\baselineskip 15pt
\topmatter
\title Sequences of analytic disks\endtitle
\author Evgeny A. Poletsky \endauthor
\thanks Partially supported by NSF grant \#DMS-9101826\endthanks
\keywords Sequences of analytic functions,  polynomial hulls, maximal functions\endkeywords
\subjclass Primary: 32C25; secondary: 32E20, 32F05 \endsubjclass
\abstract The subject considered in this paper has, at least, three points 
of interest. Suppose that we have a sequence of one-dimensional analytic 
varieties in a domain in $\Bbb C^n$. The cluster of this sequence consists 
from all points in the domains such that every neighbourhood of such points 
intersects with infinitely many different varieties. The first question is: 
what analytic properties does the cluster inherit from varieties? We give 
a sufficient criterion when the cluster contains an analytic disk, but it 
follows from examples of Stolzenberg and Wermer that, in general, clusters 
can contain no analytic disks. So we study algebras of continuous function 
on clusters, which can be approximated by holomorphic functions or 
polynomials, and show that this algebras possess some analytic properties in 
all but explicitly pathological and uninteresting cases. \par
Secondly, we apply and results about clusters to polynomial hulls and maximal 
functions, finding remnants of analytic structures there too.\par
And, finally, due to more and more frequent appearances of analytic disks 
as tools in complex analysis, it seems to be interesting to look at their
sequences to establish terminology, basic notation and properties.
\endabstract
\address Department of Mathematics,  215 Carnegie Hall, 
Syracuse University,  Syracuse, NY 13244\endaddress
\endtopmatter
\heading \S 0. Introduction\endheading
Let us consider the following example.\par
Let $A_j$ be a sequence of one-dimensional irreducible analytic subvarieties 
of 
a bounded domain $D$ in $\Bbb C^n$. {\it The cluster } of this sequence is 
the set $A$ of all points in 
$D$ such that every ball, centered at a point $z\in A$, has non-empty 
intersections  with infinitely many different varieties $A_j$. The question 
is: Does $A$ remember its analytic origin? And if it does, what kind of 
analytic geometry it inherits?\par
By the uniformization theorem for each $j$ and for each point $z\in A_j$ 
there is a holomorphic mapping $f_j$ of the unit disk $U$ into $D$ such that 
$f(0)=z$ and $f_j(U)=A_j$. Moreover, all asymptotic values of $f_j$ belong 
to the boundary $\bd D$ of $D$. It is more convenient to consider sequences 
of functions $f_j$ instead of sequences of varieties $A_j$ because we can use 
the greater power of complex analysis. Also this replacement brings into the 
consideration other important cases which we are going to discuss now.\par 
Let $X$ be a compact in $\Bbb C^n$ and $\hat X$ is its polynomial hull. If 
$z$ is a point of $\hat X$ then (see Th. 5.1 below) there is a sequence of 
holomorphic mappings $f_j:U\to\Bbb C^n$ such that $f_j(0)=z$ and boundary 
values of $f_j$ belong to an arbitrary fixed neighbourhood $V$ of $X$ 
everywhere on the unit circle $S$ except of a set of measure $\eps_j\to0$. 
Moreover, the cluster of the sequence $f_j(U)$ belongs to $\hat X$. Has 
the cluster an analytic structure in some sense?\par
The same question arises when one studies maximal function on the domain $D$ 
in the pluripotential theory. It can be proved that values of these function 
can be obtained as solutions of variational problems for some functional 
$\Hvi$ on the set of all holomorphic mappings of $U$ into $D$. If $u(z)$ is 
the extremal value of $\Hvi$ then we can find the sequence of holomorphic 
mappings $f_j$ of $U$ such that $f_j(0)=z$ and $\Hvi(f_j)\to u$. Again, as 
before, boundary values of $f_j$ lie near $\bd D$ except of the set of a 
small measure. And, again, it is interesting to know how the maximal function 
behaves on the cluster of $f_j(U)$. \par
The question about the existence of analytic structure is, of course, 
ambiguous. To make it more specific we can ask whether the cluster has an 
analytic disk inside or not? This was the first question being studied, and 
Stolzenberg \cite{12} and Wermer \cite{14} provided examples of clusters 
which do not contain any non-trivial analytic disks. These are impressive 
examples which got titles of "clusters without analytic structure". We give 
several sufficient conditions for a cluster to contain an analytic disk. But 
one of main goals of this paper is a generalization of a result of Goldmann. 
In 
\cite{6} Goldmann discovered that the Wermer's example has the uniqueness 
property: if a continuous function can be approximated by polynomials on this 
cluster and is equal to zero on some open set then it is equal to zero 
everywhere. So some analyticity stays in the cluster. We show in this paper 
that the uniqueness property holds in many quite general cases.\par
We also apply our results to maximal functions. Bedford and Kalka proved in 
\cite{1} that if a such function is smooth then through every point one can 
draw an analytic disk such that the restriction of the function to this disk 
is harmonic. For continuous functions the result does not hold in this form, 
but if we replace analytic disks by clusters then everything is true.\par
There is one more application of introduced notions. Let 
$\{f_j:U\to\Bbb C^n\}$ be a uniformly bounded sequence of holomorphic 
mappings with weakly converging  push-forward measures $(f_j)_*m=\mu(dz,f_j)$ 
of the measure $m=d\th/2\pi$ on the unit circle. The weak limit of such 
measures is, evidently, a Jensen measure. Conversely, as it was proved in 
\cite{3}, every Jensen measure is such a limit. We do not explore this 
observation here concentrating on geometric properties of clusters. 
\par
And another goal of this paper is to develop the nomenclature for this kind 
of species, establish basic properties of clusters, and to show that they 
have quite ample geometry.\par
I am very grateful to Norm Levenberg for useful discussions.
\heading \S 1. Holomorphic measures\endheading
We shall denote by $B^n(a,r) $ and $S^n(a,r) $ the ball and the sphere of 
radius $r$ with center $a$ lying in the complex $n$-dimensional space 
$\Bbb C^n$. We shall omit $a$  and $r $ when $a=0$ or $r=1$ and the index 
$n$ when this does not lead to a misunderstanding. We write $U$ for $B^1$  
and $S$ for $S^1$. A standard Lebesgue measure in $\Bbb C^n$ will be denoted 
by $m(dz)$.\par
If $E$ is a set in the closure $\ovr{D}$ a domain $D\sbs\Bbb C$ then we 
define its harmonic measure in the following way. If $E\sbs\ovr D$ is an open 
set then the {\it harmonic measure } $\om(z,E,D)$ of the set $E$ with respect 
to $D$ is equal to the infimum of all positive superharmonic functions on $D$ 
which are continuous and greater than 1 on $E$. If $F$ is an arbitrary set in 
$D$ then we define the harmonic measure $\om(z,F,D)$ to be equal to 
$\inf \om(z,E,D)$, where the infimum is taken over all open sets $E$ 
containing $F$. This function may be not upper  semicontinuous, so we take 
its regularization $$\varlimsup\limits_{w\to z}\om(z,F,D),$$ which  is 
superharmonic and will be denoted also by $\om(z,F,D)$.\par 
Let $f$ be a bounded holomorphic mapping of the unit disk into $\Bbb C^n$. 
Since $f$ has radial limit values almost everywhere on $S$, we may consider 
$f$ as a measurable mapping of the closed unit disk $\ovr U$ and define for 
a Borel set $E\sbs\Bbb C^n$ a set function 
$$\mu(E,f)=\frac1{2\pi}m(f^{-1}(E)\cap S)=\om(0,f^{-1}(E)\cap S,U).$$ This is 
a regular Borel measure on $\Bbb C^n$ and, if $f$ is not constant, 
$\mu(dz,f)$ has no atoms, i.e. the measure $\mu$ of 
each point in $\Bbb C^n$ is zero. Let $X$ be the support of $\mu$. 
Evidently, $f(U)$ belongs to the 
polynomial hull $\hat X$ of $X$, but $f(\ovr U)$ can be less than $\hat X$ as 
the example of Wermer in \cite{13} shows. The 
point $z_f=f(0)$ is uniquely determined 
as a point such that $$\int p(z)\,\mu(dz,f)=p(z_f)$$for every polynomial $p$ 
on $\Bbb C^n$. In particular, $$\int z\,\mu(dz,f)=z_f.$$We shall 
call the measure $\mu(dz,f)$ by {\it a holomorphic 
measure, } the point $z_f$ by {\it the center of the holomorphic measure } 
$\mu(dz,f)$. {\it An analytic disk, } corresponding to this holomorphic 
measure, is the image $f(U)$ of $f$.\par
Since the author could not find any references in the literature concerning 
such measures, we shall describe some of their properties.\par
If $f=(f_1,f_2,\dots,f_n)$ then we may consider measures $\mu(dz,f_j)$ on 
$\Bbb C$ and, evidently, for a set $E$, lying in the $j$-th coordinate 
complex plane, we have $\mu(E,f_j)=\mu(E\times\Bbb C^{n-1},f)$. Thus, the 
equality $\mu(dz,f)=\mu(dz,g)$ implies that $\mu(dz,f_j)=\mu(dz,g_j)$. The 
converse statement is, in general, false.\par
If analytic disk $f(U)$ is given, we may consider the mapping $f$ as a 
parameterization of $f(U)$. We have the natural method of changing parameter
by composing $f$ with a mapping $h:U\to U$. Such a reparameterization does 
not change the center of the analytic disk if $h(0)=0$ and does not change 
the holomorphic measure if and only if $h$ is inner, as the following theorem 
shows.
\proclaim{ Theorem 1.1 } For a holomorphic function $h:U\to U$ with $h(0)=0$ 
measures $\mu(dz,f)$ and $\mu(dz,f\circ h)$ are equal if and only if $h$ is 
an inner function, i.e. $|h|=1$ a.e. on $S$.\endproclaim
\demo{Proof} If $E$ is a Borel set in $\Bbb C^n$ then the measure 
$\mu(E,f)$ is equal to the  measure of the set $E_1=f^{-1}(E)\cap S$. For the 
mapping $g=f\circ h$ the measure $\mu(E,g)$ is equal to the measure of the 
set $E_2=g^{-1}(E)\cap S=h^{-1}(E_1)$. Since inner functions, mapping the 
origin into the origin, preserve measures on $S$ (see \cite{10, 19.3.4}), 
$m(E_2)=m(E_1)$.\footnote{Prof. W. Rudin told me about this reference.}\par
Conversely, for a holomorphic mapping $h$ of $U$ into $U$, we 
consider the set $A$ in $S$, where values of radial limits of $h$ belongs to 
$U$. If $m(A)=0$ then $h$ is inner, so we may assume that $m(A)>0$. Let $A_*$ 
be the image of $A$ under the mapping $g=f\circ h$. We denote by $A_1$ the 
set $f^{-1}(A_*)\cap S$ and by $A_2$ the set $g^{-1}(A_*)\cap S$. Evidently, 
$A\sbs A_2$ and we denote by $A'$ the set $A_2\sm A'$. Then 
$$m(A_1)=m(A_2)=m(A)+m(A').$$ Let $u$ be the harmonic function equal to 1 on 
$A_1$ and to zero on $S\sm A_1$ and let $v=u\circ h$. Then 
$u(0)=v(0)=m(A_k),\,k=1,2.$ But the function $v$ is equal to 0 on 
$S\sm(A\cup A')$, to 1 on $A'$ and is less then 1 on $A$. Therefore, 
$$m(A_2)=m(A)+m(A')>v(0)=m(A_2).\qed$$
\enddemo
Holomorphic measures are related to the pluripotential theory as the 
following simple reasoning shows. A set $E\in\Bbb C^n$ is the set of 
{\it universal holomorphic measure zero } if for every nonconstant 
holomorphic mapping $f:U\to\Bbb C^n$ the measure $\mu(E,f)=0$. 
Such sets exist, for example, all sets, containing only one 
point, are of universal holomorphic measure zero. In 20-th Lusin introduced 
the notion of $P$-set. A $P${\it -set } $E$ is a set in $\Bbb R^n$ such that 
$\mu(E)=0$ for every finite Borel regular measure without atoms. There are 
non-trivial examples of $P$-sets. 
\proclaim{Theorem 1.2 } \roster\item Every $P$-set in $\Bbb C^n$ is a set of 
universal holomorphic measure zero.\item Every set of universal holomorphic 
measure zero in $\Bbb C^n$ is pluripolar.
\item A Borel set $E$ in the complex plane $\Bbb C$ is a set of universal 
holomorphic measure zero if and only if it is polar.\endroster
\endproclaim
\demo{Proof} 1) is trivial.\par
Part 2) follows from Corollary 2.1.2 of \cite{9}.\par
Part 3) follows from part 2) and Theorem 2.16 of \cite{4}, claiming that 
polar sets has holomorphic measure zero for every non-constant bounded 
mapping $f:U\to\Bbb C$.\qed\enddemo
In $\Bbb C^n$ the latter statement does not hold. For example, take 
$E=\{(z,w):\,||(z,w)||<1,\,w=0\}$, which is pluri-polar, but $\mu(E,f)\ne 0$ 
for $f(\zeta)=(\zeta,0)$. It is interesting to describe all sets of universal 
holomorphic measure 0 in $\Bbb C^n$. Evidently,
\proclaim{Theorem 1.3} Let $A$ be a smooth complex analytic curve in a 
domain $D\sbs\Bbb C^n$. A set $E\sbs A$ is a set of universal holomorphic 
measure zero if and only if $E$ is polar on $A$.\endproclaim
\demo{Proof} Follows from part 3) of Theorem 1.2\qed\enddemo
Of course, we must be aware that a set $E$ can belong to supports of two 
holomorphic measures $\mu(dz,f)$ and $\mu(dz,g)$, and $\mu(E,f)\ne0$ but 
$\mu(dz,g)=0$.
\subheading{ Example 1.1 } Let $E$ be a closed set on $[-1,1]$ with linear 
measure equal to zero but non-polar as a set in $\Bbb C$. Let $f$ be a 
conformal mapping of $U$ onto the upper half of $U$ and let $g$ be the 
universal covering of $U\sm E$. Then $\mu(E,f)=0$ and $\mu(E,g)>0$.\par 
It is interesting to find out when two holomorphic have equal holomorphic 
measures. It is easy to construct an example of two holomorphic mappings 
which measures are absolutely continuous with respect to each other, but not 
equal. 
\heading \S 2. Jensen measures.\endheading
{\it A Jensen measure on $\Bbb C^n$with barycenter $z_0\in\Bbb C^n$ }  is a 
regular Borel measure $\mu$ with compact support such that 
$\mu(\Bbb C^n)=1$ and for every plurisubharmonic function $u$ on $\Bbb C^n$ 
$$u(z_0)\le\int u(z)\mu(dz).$$It is clear that every holomorphic 
measure $\mu(dz,f)$ is a Jensen measure with barycenter $z_f$.\par
Let $L=\{f_j\}$ be a sequence of uniformly bounded holomorphic mappings of 
$U$ with weakly converging holomorphic measures $\mu(dz,f_j)$, i.e. for every 
continuous function $\hvi(z)$ on $\Bbb C^n$ we have 
$$\lim\limits_{j\to\infty}\int\hvi(z)\,\mu(dz,f_j)=\int\hvi(z)\,\mu(dz)$$ 
where $\mu$ is a regular Borel measure on $\Bbb C^n$. Then centers $z_{f_j}$ 
converge to a point $z_0$ and the measure $\mu=\mu(dz,L)$ is a Jensen measure 
with barycenter $z_0$. The following theorem, claiming that every Jensen 
measure on $\Bbb C^n$ can be obtained as a weak limit of holomorphic 
measures, was proved, basically, in \cite{3}. Since we could not find the 
theorem in the form we need in \cite{3} and for the readers' convenience we 
shall give another proof of this result.
\proclaim{Theorem 2.1} Let $D$ be a domain in $\Bbb C^n$ and let $\mu$ be 
a Jensen measure with barycenter $z_0\in D$ and support compactly belonging 
to $D$. Then there is 
a sequence $L=\{f_j\}$ of holomorphic mappings of $U$ into $D$ with 
holomorphic measures $\mu(dz,f_j)$ weakly converging to $\mu$.\endproclaim 
\demo{Proof} Let $z_0=0$. Following \cite{3} we shall prove that the set 
$\bl P$ of all weak limits of sequences $\mu(dz,f_j)$, $f_j:U\to D$, 
$|f_j|<M$, and $\lim f_j(0)=0$, is convex. Let $\mu(dz,f_j)$ converges to 
$\mu$ and $\mu(dz,g_j)$ converges to $\nu$. We shall prove that 
$\al\mu+(1-\al)\nu\in\bl P$. 
We may assume that all mappings $f_j$ and $g_j$ are holomorphic in a 
neighbourhood of $\ovr U$ and $f_j(0)=g_j(0)=0$. We consider mappings 
$$h_j(\zeta)=f_j(\zeta)+g_j\left(\frac{r_j}{\zeta}\right),$$defined on 
rings $R_j=\{\zeta:\,r_j\le\zeta\le1\}$. Since either $|\zeta|$ or 
$r_j/|\zeta|$ is less than $\sqrt{r_j}$ on $R_j$, we see that for every $j$ 
and every $\eps_j>0$ the set $h(R_j)$ belongs to the $\eps_j$-neighborhood 
of $f_j(U)\cup g_j(U)$ when $r_j$ is sufficiently small. In particular, 
we can find $r_j$ such that $h_j(R_j)\sbs D$. Let $e_j(\zeta)$ be the 
composition of $e^\xi$ and conformal equivalence of the unit disk and 
the strip $\{\xi:\ln r_j<\bl{Re}\xi<0\}$, mapping the origin into 
$\xi=(1-\al)\ln r_j$ and the point $\zeta=1$ into the origin. Evidently, 
mappings $e_j$ map  the arc $\gm_1=\{e^{i\th}:-\pi\al<\th<\pi\al\}$ onto the 
unit circle $S$ and its complement $\gm_2$ on the circle $r_jS$. We consider 
mappings $p_j=h_j\circ e_j$ of $U$ into $D$.\par
If $\hvi$ is a continuous function on $\Bbb C^n$ then 
$$\int \hvi(z)\mu(dz,p_j)=
\frac{1}{2\pi}\il0^{2\pi}\hvi(h_j(e_j(e^{i\th}))+
g_j(r_je_j^{-1}(e^{i\th})))\,d\th.$$
Because of symmetry it is enough to estimate the integral 
$$\frac{1}{2\pi}\il{-\pi\al}^{\pi\al}\hvi(f_j(e_j(e^{i\th}))+
g_j(r_je_j^{-1}(e^{i\th})))\,d\th.$$
Since $|g_j(r_je_j^{-1}(e^{i\th}))|\le Mr_j$ for $\th\in\gm_1$ the previous 
integral is equal to 
$$\frac{1}{2\pi}\il{-\pi\al}^{\pi\al}\hvi(f_j(e_j(e^{i\th})))\,d\th+\dl_j,$$
where the sequence of numbers $\{\dl_j\}$ converges to 0 as $j\to\infty$. If 
$u_j$ are harmonic functions on the unit disk with boundary values 
$\hvi(f_j(e^{i\th}))$ then they are uniformly bounded and, therefore, 
$|u_j(\zeta)-u_j(0)|<K|\zeta|.$ Thus, 
$$u_j(e_j(0))=\frac{1}{2\pi}\il0^{2\pi}u_j(e_j(e^{i\th}))\,d\th=
\frac{1}{2\pi}\il{-\pi\al}^{\pi\al}\hvi(f_j(e_j(e^{i\th})))\,d\th+
\frac{1}{2\pi}\il{\gm_2}u_j(e_j(e^{i\th}))\,d\th.$$
Since $|u_j(e_j(0))-u_j(0)|<Kr_j^{1-\al}$ and 
$|u_j(e_j(e^{i\th}))-u_j(0)|<Kr_j$ for $\th\in\gm_2$, we see that 
$$\frac{1}{2\pi}\il{-\pi\al}^{\pi\al}\hvi(f_j(e_j(e^{i\th})))\,d\th=
\al u_j(0)+k_j,$$ where numbers $k_j$ converge to 0. Therefore, 
$$\int \hvi(z)\mu(dz,p_j)\to\al\int \hvi(z)\,\mu_1(dz)+
(1-\al)\int \hvi(z)\mu_2(dz),$$and we proved the convexity of $\bl P$.\par
Clearly, the set $\bl P$ is closed. If $\mu$ is a Jensen measure with 
barycenter $0\in D$ which is not in $\bl P$ 
then there is a continuous function $\hvi$ on $\Bbb C^n$ such that 
$$\int \hvi(z)\,\mu(dz)=a<b\le\int \hvi(z)\,\nu(dz)$$for all measures 
$\nu\in\bl P$. It was proved in \cite{8} (see, also, \cite{3}) that if 
$v(z)$ is a maximal plurisubharmonic function on $D$ which is less than 
$\hvi$ then the infimum over all $\nu\in\bl P$ of $$\int\hvi(z)\,\nu(dz)$$ is equal to $v(0)$.  
Since any plurisubharminic function can be approximated by continuous plurisubharmonic functions on compacts and since the support of $\mu$ is a compact in $D$ we see that $$a=\int \hvi(z)\,\mu(dz)\ge\int v(z)\,\mu(dz)\ge v(0)\ge b,$$and we 
get the contradiction.\qed\enddemo
\heading \S 3 Leaves\endheading
Let $L=\{f_j\}$ be an uniformly bounded sequence of holomorphic mappings of 
the unit disk into $\Bbb C^n$, i.e. all sets $f_j(\ovr U)$ belong to a ball 
$D=B(0,R)$ for some $R>0$. We introduce the cluster $\cl L$ of $L$, as the 
set of all points $z\in\Bbb C^n$ such that there is a sequence of points 
$\zeta_{j_k}\in U$ with $\lim\limits_{k\to\infty} f_{j_k}(z_{j_k})=z$.\par
Clusters can be quite messy and pathological. Our goal is to find inside of 
them something nice, and to make the first move at this direction we shall 
use weak compactness of  holomorphic measures.\par
By the Alouglou theorem \cite{5, V.4.2},  the set of measures $\mu(dz,f_j)$ 
is compact with respect to the weak topology defined 
by integrals of continuous functions. Therefore, every sequence $L=\{f_j\}$ 
of uniformly bounded holomorphic mappings contains a subsequence $\{f_{j_k}\}$ 
with weakly converging holomorphic measures $\mu(dz,f_{j_k})$. We limit 
ourselves to studies of such sequences.\par
A sequence of uniformly bounded holomorphic mappings $L=\{f_j\}$ is called 
{\it a leaf } if measures $\mu(dz,f_j)$ are weakly converging to the 
Jensen measure $\mu(dz)$ of $L$, which will be denoted by $\mu(dz,L)$. 
Evidently, centers $z_{f_j}$ of analytic disks $f_j$ converge to a point 
$z_L$, which we shall call {\it the center } of $L$.\par
So we see that through every point of a general cluster passes a leaf made 
from the same sequence. But leaves can 
still be quite pathological. For example, let $h_j$ be a conformal mapping 
of $U$ onto the strip $D_j=\{\zeta=x+iy\in\Bbb C:|x|<1,\,|y|<1/j\}$ and $g_j$ 
be any sequence of holomorphic mappings of $D_j$ (which can approximate any 
continuous function on [-1,1]). Then the sequence $f_j=g_j\circ h_j$ has 
as its cluster any cluster of a sequence of continuous mappings of [-1,1], 
which can be extremely unpleasant. To avoid the consideration of such 
pathologies (which is not a great loss as we show later) we should introduce 
a couple of definitions.\par
A point $z\in\cl L$ is called {\it nonessential} if there is a number $r>0$ 
such that $$\lim\limits_{j\to\infty}\om(0,f_j^{-1}(B^n(z,r)),U)=0.$$
Other points in $\cl L$ are called {\it essential. } We shall denote the set 
of essential points of $L$ by $\ess L$. Evidently, this set is closed.\par
Any leaf  $L=\{f_j\}$ determines the set $\supp L$ equal to the support of 
the measure $\mu(dz,L)$. Evidently $\supp L$  is compact, every point 
$z\in\supp L$ is essential, and $\supp L=\{z_L\}$ if and only if 
$\ess L=\{z_L\}$.\par
The following theorem claims that we can get rid of non-essential points  
without losing too much.
\proclaim{Theorem 3.1} For every uniformly bounded sequence $L=\{f_j\}$ of 
holomorphic mappings of $U$ into $\Bbb C^n$ there is another sequence of such 
mappings $M=\{g_j\}$ such that:\roster\item $\cl M$ is the set of all 
essential points of $\cl L$;\item all points of $\cl M$ are essential;
\item if $L$ is a leaf then $M$ is a leaf and $\mu(dz,L)=\mu(dz,M)$.
\endroster\endproclaim
\demo{Proof} For every nonessential point $z$ of $\cl L$ we can find a ball 
$B^n(z,r(z))$ such that
$$\lim\limits_{j\to\infty}\om(0,f_j^{-1}(B^n(z,2r(z))),U)=0.$$
Let us choose a countable covering of the set of nonessential points by 
balls $B_k=B(z_k,r_k)$, where $z_k$ is a nonessential point and 
$r_k=r(z_k)$. We denote by $C_m$ the union of first $m$ closed balls 
$\ovr B_k,\,k=1,\dots,m$. For every integer $m$ there is an integer $j(m)$ 
such that $$\om(0,f_j^{-1}(C_m),U)<\frac1m$$when $j\ge j(m)$. For 
$j(m)\le j<j(m+1)$ we consider open sets $D_{mj}=f_j^{-1}(\Bbb C^n\sm C_m)$ 
in $U$. Each of $D_{mj}$ contains the origin and let $p_j$ be an universal 
holomorphic covering mapping of $D_{mj}$ by $U$ such that $p_j(0)=0$.\par
The mapping $p_j$ has radial limits a.e. on $S$ and, by Lindel\"of theorem, 
$g_j(\zeta)=f_j(p_j(\zeta))$ for almost all points $\zeta\in S$. If $E$ is 
a subset of $U$ then, evidently, 
$$\om(0,p_j^{-1}(E),U)\ge \om(0,E,U)-\frac1m.\tag1$$
We define $g_j=f_j\circ p_j$ and let $M=\{g_j\}$. Evidently, $\cl M$ is the 
set of all essential points of $L$.\par
If $z$ is an essential point of $L$ then for every ball $B(z,r),\,r>0$ there 
is a subsequence $\{f_{j_k}\}$ such that 
$$\om(0,f^{-1}_{j_k}(B(z,r)),U)>a(r)>0.$$By (1) 
$$\om(0,g^{-1}_{j_k}(B(z,r)),U)>a(r)-\frac1m\to a(r)>0.$$
So all points of $\cl M$ are essential. If $L$ is a leaf, the same reasoning 
shows that $\mu(dz,g_j)$ weakly converge to the same measure $\mu(dz,L)$ 
and the sequence $M=\{g_j\}$ is a leaf with $\mu(dz,M)=\mu(dz,L).$\qed\enddemo
As the first application this theorem tells us that the set of essential
points is quite big.
\proclaim{Corollary 3.1} For every leaf $L$ the set $\ess L$ is connected.
\endproclaim
\demo{Proof} By Theorem 3.1 $\ess L=\cl M$ which is connected.\qed
\enddemo
A sequence $M=\{g_j\}$ is called {\it perfect } if all points of $\cl M$ are 
essential. A perfect leaf $M$ satisfying conclusions of Theorem 3.1 for 
some sequence $L$ is called {\it a perfection } of $L$. \par 
Let $V$ be an open set in $\Bbb C^n$ and $L=\{f_j\}$ be a sequence. For a 
point $z\in f_j(U)$ we define the function $\om(z,V,f_j(U))$ to be equal to 
the maximum of $\om(\zeta,f_j^{-1}(V),U)$ for all $\zeta$ such that 
$f_j(\zeta)=z$. We set $\om(z,V,f_j(U))$ to be equal to 0 if 
$z\not\in f_j(U)$. We introduce the harmonic measure on clusters as 
an upper semicontinuous function 
$$\om(z,V,L)=\varlimsup\limits_{w\to z,\,j\to\infty}\om(w,V,f_j(U))$$
on $\cl L$.\par
If $E$ is arbitrary set in $\Bbb C^n$ we let $$\om(z,E,L)=\inf \om(z,V,L),$$
where the infimum is taken over all open sets $V$ containing $E$.\par
\proclaim{Lemma 3.1} If a point $z\in\cl L$ is essential then for  
every open set $E\sbs\cl L$, containing $z$, $\om(z_L,E,L)>0$.
\endproclaim
\demo{Proof} Let $E$ be an open set in $\cl L$. We can find $r>0$ such 
that the intersection of the ball $B_2=B(z,2r)$ and $\cl L$ belongs to $E$. 
For an arbitrary $\dl>0$ there is an open set $V\sbs\Bbb C^n$ such that 
$V\cap\cl L=E$ and $\om(z_L,V,L)<\om(z_L,E,L)+\dl$. There is an integer $k$ 
such that $f_j(\zeta)\in V$ if $f_j(\zeta)\in B(z,r)$ and $j>k$. Therefore, 
$\om(z_L,V,L)>\om(z_L,B(z,r),L)$. Since $z$ is essential, 
$\om(z_L,B(z,r),L)>0$.\qed\enddemo 
\proclaim{Theorem 3.2} Let $L=\{f_j\}$ be a perfect leaf and $X=\supp L$. 
Then the harmonic measure $\om(z,X,L)\equiv1$ 
on $\cl L$.\endproclaim
\demo{Proof} Let $V$ be an open neighbourhood of $X$. We must prove that for 
every $z\in\cl L$ and for every neighbourhood $W$ of $z$ there is a sequence 
of points $a_{j_k}\in U$ such that:\roster\item $f_{j_k}(a_{j_k})\in W$;\item 
$$\lim\limits_{k\to\infty}\om(a_{j_k},f^{-1}_{j_k}(V),U)=1.$$\endroster
This is evident if $z_L\in W$ because $\om(z_L,V,L)=\mu(\Bbb C^n,L)=1$. So we 
may assume that $z_L\not\in W$.\par
Let $W_j=f_j^{-1}(W)$ and let $B_1\Sbs W$ be a ball centered in $z$. Since 
$z$ is essential there is a subsequence $j_k$ such that  for sets 
$Z_k=f^{-1}_{j_k}(\ovr B_1)$ the harmonic measure $$\om(0,Z_k,U)>\dl>0.$$
Let $V_k=f^{-1}_{j_k}(V)\cap S$ and $v_k(\zeta)=\om(\zeta,V_k,U)$. Suppose 
that $v_k$ is less than $\eps<1$ on $W_{j_k}$. If $g_k$ is an universal cover 
of $U\sm Z_k$ such that $g_k(0)=0$, then radial limits of $g_k$ are in $Z_k$ 
on a set of measure greater than $\dl$. Therefore, for the function 
$u_k=v_k\circ g_k$ we have 
$$u_k(0)<(1-\dl)+\eps\dl=1-\dl(1-\eps)<1.$$
But the sequence $u_k(0)=v_k(0)$ converges to 1, so there are an integer 
$k$ and a point $a_{j_k}\in W_{j_k}$ such that 
$$\om(a_{j_k},f^{-1}_{j_k}(V),U)>\eps.$$ Since $\eps<1$ is arbitrary we get 
our theorem.\qed
\enddemo
A point $z\in\cl L$ is called {\it totally essential } if for every ball 
$B_1=B(z,r)$ the lower limit of the sequence $\om(0,f^{-1}_j(B_1),U)$ is 
greater than 0. Evidently, the set of totally essential points is closed and 
all points in $\supp L$ are totally essential. A leaf is {\it totally 
perfect } if all point of its cluster are totally essential.
\proclaim{Theorem 3.3} If $L=\{f_j\}$ is a perfect leaf then for every 
$z\in\cl L$ 
there is a totally perfect leaf $M$ such that $z_M=z_L$, $z\in\cl L$, 
$\mu(dz,M)=\mu(dz,L)$, and $\cl M$ belongs to $\cl L$ and contains all 
totally essential points of $L$.\endproclaim
\demo{Proof} Since $z$ is an essential point of $L$, it is easy to see that 
there is a subsequence $\{f_{j_k}\}$ such that $z$ is totally essential for 
this subsequence. Its perfection contains $z$ and the set $T$ of all totally 
essential points of $L$, which will continue to be totally essential. So 
we may assume that $z$ is a totally essential point of $L$.\par
Let $n_1>1$ be the minimal positive integer such that the set 
$$E_1=\{z\in\cl L:\,\opr{dist}(z,T)\ge\frac1{n_1}\}$$is non-empty. Since 
$E_1$ is compact there is a finite cover of $E_1$ by balls $B_{1m}$ of radius 
less $1/2n_1$, centered at points of $E_1$ and such that for each ball 
$B_{1m}$ there is a subsequence $\{f_{j_k}\}$ with 
$$\lim\limits_{k\to\infty}\om(0,f^{-1}_{j_k}(B_{1m}),U)=0.$$
We take such a subsequence $\{f_{1j}\}$ for the ball $B_{11}$. For its 
perfection $L_1$ the set $T_1$ of all totally essential points contains $T$ 
and the $\cl L_1$ does not intersect $B_{11}$. If the ball $B_{12}$ contains 
points of $T_1$ we leave it at peace, otherwise we repeat the procedure for 
the sequence $\{f_{1j}\}$. After finite number of steps we exhaust all balls 
$B_{1m}$ and come to the leaf $L_l$ with the set $T_l$ of totally essential 
points. Every point of $\cl L_l$ lies within the distance $1/n_1$ from points 
of $T_l$. So if $n_2$ is the minimal positive integer such that the set 
$$E_2=\{z\in\cl L_l:\,\opr{dist}(z,T_l)>\frac1{n_2}\}$$is non-empty, then 
$n_2>n_1$. We again take a finite cover by balls of radius less than $1/2n_2$ 
and repeat the previous construction.\par
At the end of the process we get a sequence of leaves $M_k=\{f_{kj}\}$ with 
perfections $L_k$ and the sequence of balls $B_k$ such that $M_1$ is a 
subsequence of $L$, $M_{k+1}$ is a subsequence $M_k$, the ball $B_k$ is 
non-essential for $M_k$.  Sets $T_k$ are increasing and their union is dense 
in the intersection of clusters of $L_k$. We take the sequence 
$\{g_j=f_{jj}\}$ and let $M$ be its perfection.\par
First of all, $z_M=z_L$ and $\mu(dz,M)=\mu(dz,L)$. Secondly, $\cl M$ belongs 
to the intersection of clusters of $L_k$ and, therefore, doesn't contain the 
union of balls $B_k$. And if the point $z\in T_k$ for some $k$ then it is 
totally essential for $\{g_j\}$. Therefore, $\cl M$ consists from totally 
essential points only.\qed\enddemo
Let us give an example of a perfect but not a totally perfect leaf.
\subheading{Example 3.1} Let $f_j=(\zeta,\zeta^j)$ be a sequence of mappings 
of $U$ into $U^2$. It is easy to see that measures $\mu(dz,f_j)$ converge to 
the measure $$\mu=\frac1{2\pi}d\al\,d\beta$$ on $S\times S$. Therefore, 
the sequence of $f_j$ determines the leaf $L=\{f_j\}$ with $\mu(dz,L)=\mu$. 
So $\supp L=S\times S$, and 
$$\cl L=\{(z_1,z_2):|z_1|<1,\,z_2=0 \text{ or } |z_1|=1,\,
|z_2|\le 1\}.$$ The leaf $L$ is totally perfect.\par
We can consider another sequence of holomorphic mappings $g_j$ equal to 
$f_j$ if $j$ is even, and to $(\zeta^j,\zeta)$ if $j$ is odd. Then 
holomorphic measures $\mu(dz,g_j)$ also converge to $\mu$, so the sequence 
$M=\{g_j\}$ is a leaf, and this leaf is perfect but not totally perfect. 
Its cluster is the union of two perfect leaves.\par
\heading \S 4. Geometry and analytic structures on leaves\endheading
If $L=\{f_j\}$ is a leaf in $\Bbb C^n$ then it is interesting to find out 
what does $\cl L$ remember about its holomorphic origin. There is no sense 
to look at leaves with trivial Jensen measures: their clusters can be as 
pathological as clusters of continuous mappings. So we assume that 
$\supp L\ne\{z_L\}$. In this case, the set $\ess L\ne\{z_L\}$ and, by 
Theorem 3.1, we can find a perfect leaf with the same support and the same 
Jensen measure, which cluster is $\ess L$. Since remaining points again 
can form extremely pathological set we shall concentrate on studies of 
perfect leaves.\par
The first question which naturally arises is the question about the existence 
of analytic disks in clusters. The negative answer to this question was 
given by Stolzenberg in \cite{12} and Wermer in \cite{14}. There are perfect 
leaves such that any holomorphic mapping of the unit disk into the cluster of 
such leaves is constant. 
Before we go further, we give the brief description of these examples. 
\subheading{Examples of Stolzenberg and Wermer} Both examples deal with 
complex analytic 1-dimensional irreducible varieties $V_j$ (Stolzenberg) and 
$W_j$ (Wermer) in the unit bidisk in $\Bbb C^2$. Each of $V_j$ or $W_j$ 
contains the origin. So, according to our approach, these sequences of 
varieties or their subsequence are leaves with supports lying on the boundary 
of $U^2$. It was proved by Stolzenberg and Wermer that clusters of these 
leaves does not contain an analytic disk.\par 
There is a difference in this examples. The example of Stolzenberg is a 
"Swiss cheese" lifted to $\Bbb C^2$, the projections of its cluster on 
both coordinate planes are nowhere dense. But if one makes a 
"Swiss cheese" on a plane then it may happen that its cluster (what is 
left after we drag out holes) consists from non-essential points only. By 
Theorem 3.1 Stolzenberg's example has a connected set of essential points 
joining the origin with the boundary of $U^2$.\par
Varieties $W_j$ in Wermer's example have proper projections on one of 
coordinate lines, say, $\{w=0\}$. Moreover, there is a countable dense 
set of points $\{p_k\}$ on this line such that for each $j$ the variety $W_j$ 
has finite number of points over each of $p_k$, and there is 
only one point in $W_j$ over the origin. This property was used by Goldmann 
in \cite{6} to show that the cluster of Wermer's example is perfect. In 
the same paper Goldmann discovered that this cluster has some analytic 
properties, which we discuss later. This discovery gave us the starting push 
for our studies. But before we move to this questions we want to give a 
sufficient conditions for the existence of analytic disks in clusters.
\proclaim{Theorem 4.1} Let $L=\{f_j\}$ be a leaf in $\Bbb C^n$. Suppose that 
there is a holomorphic function $h$, defined on a neighbourhood $V$ of 
$\cl V$, and numbers $b<1$ and $\al>0$ such that:\roster\item $|h|$ and 
$||\nabla h||$ are less than $M<\infty$ on $V$; 
\item $|h(z_L)|\le b<1$;\item for every integer $j\ge1$ there is an arc 
$\gm_j$ with the length greater than $\al$ such that 
$$\varliminf\limits_{r\to1}|h(f_j(r\zeta))|\ge1$$for all points $\zeta$ in 
$\gm_j$.\endroster
Then there is a non-constant holomorphic mapping 
$g=(g_1,\dots,g_n):U\to\cl L$ such that at least one of functions 
$g_k,\,1\le k\le n$, covers univalently a disk of a positive radius 
$r=r(M,b,\al)$ in $\Bbb C$.\endproclaim 
\demo{Proof} Let $S(c)$ be a set of holomorphic functions $u$ on $U$ with the 
Bloch norm $$||u||_B=\sup\limits_{\zeta\in U}|u'(\zeta)|(1-|\zeta|^2)\ge c.$$
By  Theorems 2 and 1 from \cite{7} functions $h_j=h\circ f_j$ belong to 
$S(c)$ for some $c=c(M,b,\al)>0$. Since $||\nabla h||<M$ there is an integer 
$k(j)$ such that $$||f_{jk(j)}||_B\ge\frac c{nM},$$where $f_{jk}$ is the 
$k$-th coordinate function of $f_j$. So for some $k$ between 1 and $n$ we can 
find an infinite subsequence $f_{j_l}$ such that:\roster\item 
$$|f'_{j_lk}|(1-|\zeta_{j_l}|^2)\ge\frac c{2nM}$$ at some point $\zeta_{j_l}$ 
in $U$;\item the sequence $f_{j_l}(\zeta_{j_l})$ converges to $z_0$.
\endroster 
Mappings 
$$f_{j_l}\left(\frac{\zeta-\zeta_{j_l}}{1-\ovr\zeta_{j_l}\zeta}\right)$$ 
converge to a non-constant mapping $g$ which image lies in $\cl U$ and 
$|g'_k(0)|\ge c/2nM$. So, by 
the Bloch theorem, $g_k$ covers univalently a disk of a positive radius 
$r=r(M,b,\al)$.\qed\enddemo
Examples of Stolzenberg and Wermer (see, also, \cite{7}) show that it is 
hard to await something really good at this situation but in the 
one-dimensional case we can do a little bit better.
\proclaim{Theorem 4.2} Let $L=\{f_j\}$ be a perfect leaf in $\Bbb C$ and 
$\Dl$ be the boundary of $\cl L$. Then $\Dl\sbs\supp L$.\endproclaim
\demo{Proof} Suppose that there is a disk $B_1=B(z_0,r)$ such that 
$z_0\in\Dl$, $r>0$ and $\mu(B_1),L)=0$. We can find a point $z$ such that 
a closed disk $\ovr B_2=\ovr B(z,s)\sbs B(z_0,r/4)$ does not intersect 
$\cl L$. Therefore, there is an integer $N$ such that 
$f_j(U)\cap B_2=\emptyset$ when $j>N$. Let $B_3=B(z,r/2)$ and let $u(z)$ be a 
harmonic function on $\Bbb C\sm B(z,s)$, equal to 1 on $\bd B_2$ and 
to 0 on $\bd B_3$. Functions $u_j(\zeta)=u(f_j(\zeta))$ will be harmonic on 
$U$. If $E_j=\{\zeta\in S: u_j(\zeta)>0\}$ then the 
sequence of $a_j=m(E_j)$ tends to 0. If $B_3=B(z_0,r/4)$ then $u(z)>c>0$ on 
$B_3$. Harmonic functions $v_j$ on $U$ with boundary values equal to 1 on 
$E_j$ and to 0 on $S\sm E_j$ are greater than $u_j$. The harmonic measure of 
the set $F_j=\{\zeta\in U: v_j(\zeta)>c\}$ is less than $a_j/c$ and, 
therefore, the harmonic measure of sets $G_j=\{\zeta\in U: u_j(\zeta)>c\}$ is 
less than $a_j/c$ which means that the point $z_0$ is nonessential.
\qed\enddemo
In general, it may happen that $\mu(\Dl\cap B(z,r),L)=0$ for some $z\in\Dl$. 
\proclaim{Corollary 4.1} Let $L=\{f_j\}$ be a leaf in $\Bbb C$ such that 
$z_L=0$  and $\supp L$ lies outside the disk $B_1=B(0,r)$ for some $r>0$. 
Then $B_1$ belongs to $\cl L$.\endproclaim
\demo{Proof} We may assume that the leaf $L$ is perfect. Since $\supp L$ does 
not intersect $B_1$, the boundary of $\cl L$ does not intersect $B_1$. So 
$B_1$ belongs to $\cl L$.\qed\enddemo
\subheading{Example 4.1} Let $f_j=(g_j,h_j)$ be uniformizations of curves 
$V_j$ from the Stolzenberg's example by analytic disks such that 
$f_j(0)=(0,0)$. We may assume that $\{f_j\}$ is a leaf $L$. Then 
$\supp L\sbs\bd U^2$ and, therefore, either the leaf $G=\{g_j\}$ or the leaf 
$H=\{h_j\}$ is non-trivial. Their perfections are also non-trivial. 
Nevertheless, their clusters do not contain any disk in $\Bbb C$. In this 
example, $\cl L=\supp L=\Dl$.\par
So the complex geometry is not inherited by leaves explicitly. To go further 
we introduce for each leaf $L$ algebras $\Cal P(L)$ and $\Cal H(L)$ 
of functions on $\cl L$, which can be approximated on $\cl L$ either by 
polynomials or by functions holomorphic in a neighbourhood of $\cl L$. Such 
algebras are called {\it analytic } if a function from the algebra is equal 
to 0 on $\cl L$ provided it is equal to 0 on some open set $V$ in $\cl L$.
\par
Our first step is to prove the two constant theorem for functions in 
$\Cal H(L)$.
\proclaim{ Theorem 4.3} Let $L=\{f_j\}$ be a sequence of uniformly bounded 
analytic disks and let $h\in\Cal H(L)$. If $|h|$ is less than $M$ on $\cl L$ 
and less than $m<M$ on an a set $E\in\cl L$, then 
$$|h(z)|<m^dM^{1-d},$$where $d=\om(z,E,L)$.\endproclaim
\demo{Proof} Let us fix some $\eps>0$. The function $h$ is the limit of 
functions $h_k$ holomorphic in neighbourhoods $D_k$ of $\cl L$. We may find 
an integer $k$ such 
that $|h_k|<M'=(1+\eps)M$ on $D_k$, $|h_k|<m'=(1+\eps)m$ on $E$, and 
$|h_k(z)-h(z)|<\eps$.\par
Let $V$ be the open set of points in $D_k$, where $|h_k|$ is less than $m'$. 
Evidently, $\om(z,V,L)\ge d$. If $V_j=f_j^{-1}(V)$ then, by the definition of 
$\om(z,V,L)$, there is a sequence of points $a_{j_l}\in U$ such that points 
$f_{j_l}(a_{j_l})$ converge to $z$ and the limit $d'$ of the sequence 
$d_{j_l}=\om(a_{j_l}, V_{j_l}, U)$ is greater than $d$. We may assume that 
$j_l=j$.\par
If $v_j=h_k\circ f_j$ then by the two constants theorem 
$$|v_j(a_j)|<(m')^{d_j}(M')^{1-d_j}$$ and 
$$|h_k(z)|<(m')^{d'}(M')^{1-d'}\le (m')^d(M')^{1-d}.$$
Therefore, $|h(z)|<m^dM^{1-d}$.\qed\enddemo
\proclaim{ Corollary 4.2 } Let $L=\{f_j\}$ be a perfect leaf in 
$\Bbb C^n$ 
and let $z_L$ be the center of $L$. If a function $h\in\Cal H(L)$ is equal 
to zero on an open set $E\sbs\cl L$ then $h(z_L)=0$.\endproclaim
\demo{Proof} This statement follows from the previous theorem because all 
points in $\cl L$ are essential and, therefore, $\om(z_l,E,L)>0$.\qed\enddemo
\proclaim{Corollary 4.3} Let $L=\{f_j\}$ be a perfect leaf in $\Bbb C^n$. 
Then: \roster\item every function $h\in\Cal H(L)$ attains its maximum modulus 
value on $\supp L$;\item cluster of $L$ belongs to the polynomial hull of 
$\supp L$.\endroster\endproclaim
\demo{Proof} Part 1) follows immediately from Theorem 4.3 and 3.2 and part 2) 
is the consequence of part 1).\qed\enddemo 
The maximum modulus principle proved just now can be obtained from the fact 
that the measure $\mu(dz,L)$ is an Jensen measure for $\Cal H(L)$ at $z_L$. 
We did not explore this approach because we intended to be as geometric as 
possible. Connections of this subject with the theory of uniform algebras are 
evident. Nevertheless, it seems to us that this theory does not allow to 
look closely at geometric structures generated by leaves. Moreover, we shall 
prove in \S 5 that polynomial hulls can be foliated by perfect leaves, so 
our method can be applied be applied to uniform algebras as well.
Let us consider the following example.
\subheading{Example 4.2} Let $L$ be a perfect leaf from Example 3.1. The 
function $h(z_1,z_2)=w$ is equal to 0 on a neighbourhood of origin and is not 
0 everywhere. So perfect leaves lack of the propagation of zeros. To 
understand this phenomenon we must look closer at the geometry of leaves.\par 
It happens that leaves generate complex substructures in a natural way which 
we are going to describe now. Let $$G_a(\zeta)=\frac{\zeta+a}{1+\ovr a\zeta}$$
be a conformal automorphism 
of the unit disk. Let $L=\{f_j\}$ be a perfect leaf. For a point  
$z\in \cl L$ we consider points $a_k\in U$ such that the a subsequence 
$f_{j_k}(a_k)$ converges to $z$. From the sequence of mappings 
$g_k=f_{j_k}\circ G_{a_k}$ we can choose a subsequence $h_m$ with weakly 
converging measures $\mu(dz,h_m)$. {\it A subleaf } $L_z$ is the perfection 
of the leaf $M=\{h_k\}$.\par  
If $M$ is a subleaf of $L$ then $\cl M\sbs\cl L$ but it may happen that 
$\supp M$ is not a subset of $\supp L$. For example, if $f_j$ are an 
universal holomorphic coverings of $U\sm\ovr B(1/2,1/j)$ by $U$ such that 
$f_j(0)=0$, then there are subleaves $M=L_{\frac12}$ such that 
$\mu(\frac12,M)=a$, where $a$ is any number between 0 and 1. A subleaf $L_z$ 
such that $\supp L_z\sbs\supp L$ is called {\it a principal subleaf.}\par
The set of all subleaves of the leaf $L$ is partially ordered by the 
inclusion relation of their clusters. A leaf $L$ is called {\it minimal at a 
point } $z\in\cl L$ if it contains no proper non-trivial subleaves passing 
through $z$. Evidently, a leaf $L$, minimal at $z_L$, is totally perfect. A 
leaf $L$ is called {\it minimal } if it is minimal at all points of 
$\cl L$. A point $z\in\cl L$ is called {\it inner } if there is a non-trivial 
subleaf $L_z$. We shall say that a sequence of leaves $\{L_j\}$ {\it 
converges } to a leaf $L$ if every neighbourhood of every point of $\cl L$ 
contains points from clusters of infinitely many different leaves $L_j$ and 
$\mu(dz,L_z)=\lim\mu(dz,L_j)$\par
In Example 4.2 there is only one subleaf passing through the origin, so $L$ 
is minimal at $z_L$. The leaf $M$ from Example 3.1 has two subleaves passing 
through the origin so it is not minimal at $z_L$. Both leaves have subleaves 
at the boundary of the bidisk so they are not minimal. These subleaves are 
analytic disks foliating the boundary. If we take the leaf $L$ and a point 
$z_0=(z_1,z_2)$ such that $|z_1|=1$ and $z_2\ne0$, $|z_2|<1$, then for each 
$j$ we can find a point $a_j\in U$ such that $f_j(a_j)$ converges to $z_0$ 
 and let $g_j=f_j\circ G_{a_j}$. We may assume that the sequence $g_j$ ia a 
leaf $L'=\{g_j\}$, otherwise we take a converging subsequence. Evidently, 
$\cl L=\cl L'$, but the leaf $L'$ is not perfect. The origin is a 
non-essential point of $M$. We shall show that the property for a leaf to be 
minimal and 
the analyticity of the algebra $\Cal H(L)$ are closely related.\par
\proclaim{Lemma 4.1}\roster\item If points $z_j\in\cl L$ 
converge to a point $z$ and measures $\mu(dz,L_{z_j})$ weakly converge to a 
measure $\mu$ then subleaves $L_{z_j}$ converge to a subleaf $L_z$;
\item if $J$ is a totally ordered set and $\Cal L=\{L_{z_j},\,j\in J\}$ is a 
set of subleaves such that $\cl L_{z_i}\sbs\cl L_{z_j}$ for $i>j$, then there 
is a countable sequence $\{L_j\}$ in $\Cal L$ converging to a subleaf which 
belongs to all subleaves in $\Cal L$;
\item if $z$ belongs to the cluster of a perfect leaf $L$ but doesn't belong 
to its support then there is a principal subleaf $L_z$ passing through $z$.
\endroster\endproclaim
\demo{Proof} 1) If leaves $L_{z_j}=\{g_{jk}=f_j\circ G_{a_{jk}}\}$ then we 
can find a subsequence $a_m=a_{j_mk_m}$ in $U$ such that measures 
$\mu(dz,g_{j_mk_m})$ converge to the measure $\mu$. Then the sequence 
$f_{j_m}\circ G_{a_m}$ defines a leaf $L_z$.\par
2) Let $G$ be the intersection of clusters of all leaves in $\Cal L$. There 
is a countable sequence $j_k$ such that the set $G$ is the intersection of 
$\cl L_{j_k}$. Let us select another subsequence with converging holomorphic 
measures. This subsequence of subleaves converges to a subleaf $L_z$ due to 
the first part of the lemma.\par
3) Suppose that there is a point $z\in\cl L$ such that $\dist(z,\supp L)=r>0$ 
and for every $L_z$ the support of $L_z$ doesn't belong to $\supp L$. This 
means that there are  numbers $a>0$ and $d>0$, and a ball $B_1$, centered 
at $z$, such that $$\lim\limits_{j\to\infty}\mu(W,f_j\circ G_\zeta)>a$$for 
every point $w$ in $B_1$ and every $\zeta\in f^{-1}(w)$ and for the set $W$ 
of points in $\Bbb C^n$, with distances to $\supp L$ greater or equal to 
$d$. Otherwise, we can find a subsequence $f_{j_k}=g_k$, a sequence 
$w_k\in\Bbb C^n$ converging to $z$, and a sequence $\zeta_k\in g^{-1}_k(w_k)$ 
such that for the complement $V$ of $\supp L$ measures 
$\mu(V,g_k\circ G_{\zeta_k})$ converge to 0, and, therefore, there is a 
subleaf with the support lying in the support of $L$.\par
Let $E_j$ be a subset of $S$ such that $f_j(\zeta)\in W$ for every 
$\zeta\in E_j$. If $D_j=f_j^{-1}(B_1)$ then the function 
$\om(\zeta,E_j,U)$ is greater than $a$ on $D_j$ when $j$ is big enough. Since 
$z$ is an essential point of $L$ then $\om(0,D_j,U)>\eps>0$ for infinitely 
many $j$. Therefore, $\om(0,E_j,U)>\eps a$ for such $j$, which means that 
$\mu(W,L)>\eps a>0$.
\qed\enddemo
Inner points on leaves play the role close to the role of inner points in  
domains. For example, if $D$ is a domain in $\Bbb C$ and $f$ is a mapping 
of $U$ onto $D$ then for the leaf $L=\{f\}$ the inner points are points in 
the interior of $D$. The third part of Lemma 4.1 also shows that every point 
of $\cl L$, which does not belong to the $\supp L$, is inner. In particular, 
for examples of Stolzenberg and Wermer all points of clusters lying in $U^2$ 
are inner.
\proclaim{Corollary 4.4} If $z$ is a point of a perfect leaf $L$ which does not belong to $\supp L$ 
then there is a minimal non-trivial leaf $L_z$.\endproclaim
\demo{Proof} We consider the set of all principal subleaves passing 
through $z$. By Lemma 4.1 every totally ordered subset of this set has a 
minimal element. So, by Zorn lemma, there is a minimal element in this set.
\qed\enddemo
\proclaim{Lemma 4.2} Let $f:U\to\Bbb C^n$ be a holomorphic mapping such that 
$||f||<R$, $f(0)=z_0$, $b=||z_0||\ne0$, and the harmonic measure 
$\om(z_0,B_1,f(U))=a>0$ for a ball $B_1=B(0,r),\,r<b$. Then for every $k>1$ 
such that $kr<b$ there are a number $c=c(k,r,R,a,b)>0$, a number 
$m=m(k,r,R,a,b),\,k>m>1$, and a point $z_1\in S(0,mr)$ where 
$\om(z_1,B_1,f(U))>c$ and $\om(z_1,\Bbb C^n\sm B(0,kr), f(U))>c$.\endproclaim 
\demo{Proof} Let $$u(z)=\frac{\ln\frac{||z||}{R}}{\ln\frac{r}{R}}.$$The 
function $u$ is plurisuperharmonic, positive on $B(0,R)$, and greater than 1 
on $B_1$. If $1<m<k$ and $d=\om(0,f^{-1}(B(0,mr),U)$ then 
$$u(z_0)\ge d\frac{\ln\frac{mr}{R}}{\ln\frac{r}{R}}$$or
$$d\le\frac{\ln\frac{b}{R}}{\ln\frac{mr}{R}}.$$
Let $v(\zeta)=\om(\zeta,f^{-1}(B_1),U)$ and let $\gm=f^{-1}(S(0,mr))$. If 
$q$ is the maximum of $v$ on $\gm$ then $qd\ge a$ or 
$$q\ge\frac ad=a\frac{\ln\frac{mr}{R}}{\ln\frac{b}{R}}\ge 
a\frac{\ln\frac{kr}{R}}{\ln\frac{b}{R}}=2s.$$
So there is a point $\zeta_1$ on $\gm$ where $v$ is greater than $s$. Let 
$z_1=f(\zeta_1)$.\par
Let $B_2=B(0,kr)$ and $p=\om(\zeta_1,f^{-1}(\Bbb C^n\sm B_2),U)$. Then 
$$u(z_1)=\frac{\ln\frac{mr}{R}}{\ln\frac{r}{R}}\ge 
s+(1-s-p)\frac{\ln\frac{kr}{R}}{\ln\frac{r}{R}}$$or
$$\ln\frac{mr}{R}\le 
s\ln\frac{r}{R}+(1-s)\ln\frac{kr}{R}-p\ln\frac{kr}{R}=(1-s)\ln k+
\ln\frac{r}{R}-p\ln\frac{kr}{R}.$$Therefore, 
$$p\ge\frac{(1-s)\ln k-\ln m}{\ln\frac{kr}{R}}.$$
If $$m=k^{(1-s)/2}$$ then $$p\ge\frac{(1-s)\ln k}{2\ln\frac{R}{kr}}=t.$$
If we take $c=\min(s,t)$ we get the proof of the lemma.\qed\enddemo
\proclaim{Theorem 4.4} For a perfect non-trivial leaf $L=\{f_j\}$:\roster
\item the set of inner points of $L$ is dense in $\cl L$;\item if a set $V$ 
is open in $\cl L$ and $z\in\cl L$ then $\om(z,W,L)=0$ for every set 
$W\Sbs V$ if and only if any non-trivial subleaf $L_z$ does not intersect $V$.
\endroster\endproclaim
\demo{Proof} 1) We assume that $||f_j||<R/2$. Since $L$ is non-trivial we can 
find a point $z_1\in\cl L$ such that $||z_1-z_L||=b>0$. For a ball 
$B_1=B(z_1,r),\,r>0$, we consider the ball $B_2=B(z_1,2r)$. Since the leaf 
$L$ is perfect we can find a subsequence $f_{j_k}$ such that 
$\om(z_L,B_1,f_{j_k}(U))>a>0$. By Lemma 4.2 there are points 
$\zeta_{j_k}\in U$ such that $f_{j_k}(\zeta_{j_k})=z_k\in S(z_1,mr),\,1<m<2$, 
and $\om(z_k,\Bbb C^n\sm B_2,f_{j_k}(U))>c>0$. From the sequence 
$$g_k=f_{j_k}\circ G_{\zeta_{j_k}}$$we can select a subsequence forming a 
leaf $N=\{h_l\}$and let $M$ be the perfection  of $N$. Then 
$z_M\in S(z_1,mr)$ and there are points of $\cl M$ lying outside $B_2$. 
Otherwise all this points are non-essential for the sequence $h_l$, so we 
can cover the compact intersection of $\Bbb C^n\sm B_2$ and $\cl N$ by a 
finite number of balls such that for each of them its harmonic measure at 
$z_M$ tends to 0 as $l$ goes to $\infty$. This means that 
$\om(z_M,\Bbb C^n\sm B_2,f_{h_l}(U))>c>0$ tends to 0 and this is the 
contradiction. So the leaf $M$ is non-trivial and, since $r$ is arbitrary, we 
get the proof of the first part of the theorem.\qed\par
2) If $\om(z,W,L)=0$ then for every sequence 
$g_k=f_{j_k}\circ G_{\zeta_{j_k}}$ 
such points $f_{j_k}(\zeta_{j_k})$ converge to $z$, all points in $V$ are 
non-essential and, hence, cannot belong to the perfection of $\{g_k\}$. So 
every subleaf $L_z$ does not intersect $W$ and, therefore, $V$.\par
Conversely, if $\om(z,W,L)>0$ for some $W\Sbs V$ then there are essential 
points in $\ovr W$ for the sequence $\{g_k\}$. Therefore, there are points of 
$V$ in the perfection of this sequence.\qed
\enddemo
\proclaim{Corollary 4.5} If a perfect leaf $L$ is minimal at $z_L$ and 
$z_L$ does not belong to $\supp L$ then 
for every open set $V$ in $\cl L$ there is an open set $W$ containing $z_L$ 
such that every function $h\in\Cal H(L)$ equal to 0 on $V$ is equal to 0 
on $W$.\endproclaim
\demo{Proof} Suppose that there is a sequence of points $z_k$ converging to 
$z_L$ such that $h(z_k)\ne0$. By Theorem 4.4 we may assume that points $z_k$ 
are inner and, therefore, there are principal subleaves 
$L_k=L_{z_k}=\{f_{kj}\}$, 
which, by the second part of Theorem 4.4, don't intersect $V$. If we make 
the perfect leaf $M$ from the sequence $\{f_{kk}\}$, then $z_M=z_L$ and $M$ 
is non-trivial because $\supp M$ belongs to $\supp L$ and, hence, doesn't 
contain $z_M$. But $\cl M$ is a proper subset because it doesn't contain $V$ 
which contradicts to the minimality of $L$ in $z_L$.\qed\enddemo 
\proclaim{Theorem 4.5} A leaf $L$ is minimal at any point $z\in D$ if and 
only if $\om(z,V,L)>0$ for every inner point $z\in\cl L$ and every open set 
$V\sbs\cl L$.
\endproclaim 
\demo{Proof} If $\om(z,V,L)=0$ for some open $V$ and some inner point $z$ 
then, by Theorem 4.4, a non-trivial subleaf $L_z$ doesn't intersect $V$ and, 
hence, $L$ is not minimal at $z$.\par
Conversely, if $L$ is not minimal at $z$ then there is an open set $V$ which 
doesn't intersect some $L_z$. By the second part of Theorem 4.4 this implies 
that $\om(z,V,L)=0$.\qed\enddemo
\proclaim{Corollary 4.6} If a leaf is minimal then the algebra $\Cal H(L)$ is 
analytic.\endproclaim
\demo{Proof} Follows immediately from Theorem 4.5, 4.4 1), and 4.3.\qed
\enddemo
In his paper \cite{6} Goldmann proves, in different terminology, that the 
leaf corresponding to the Wermer's example is minimal. So such leaves exist. 
Nevertheless, it seems that in general situation minimality does not happen 
frequently. If we consider the leaf $L$ from Example 3.1 then the only part 
of this leaf where we can guarantee that the algebra $\Cal H(L)$ is analytic 
is the disk $\{z_2=0\}$. To find an approach to this situation  let us 
introduce the new definition. If $L=\{f_j\}$ is a leaf then {\it the midrib } 
of $L$ is the closure of the set of all points $z$ in $\cl L$ such that 
$\om(z,W,L)>0$ for every neighbourhood $W$ of $z_L$. For the leaf $L$ 
mentioned above the midrib is the disk $\{z_2=0\}$.
\proclaim{ Corollary 4.7} If a perfect leaf $L$ is minimal at $z_L$ and $z_L$ 
does not belongs to $\supp L$ then every holomorphic function $h$ on $L$ 
equal to 0 on an open set $V$ in $\cl L$ is equal to 0 on the midrib of $L$.
\endproclaim
\demo{Proof} By Corollary 4.5 the function $h$ is equal to 0 a neighbourhood 
$W$ of $z_L$ and, by two constants theorem, is equal to 0 on the midrib.\qed
\enddemo
I do not know whether midribs are leaves or not, and I also do not know 
whether they are always nontrivial. The following theorem, which can be 
applied to the Wermer's example, gives sufficient criterion for the existence 
of a nontrivial midrib.
\proclaim{Theorem 4.6} Let $L=\{f_j\}$ be a leaf and let $h$ be a function 
from $\Cal H(L)$ such that $h(z_L)\ne h(z)$ when $z\ne z_L$ is in $\cl L$. We 
suppose that $z_L$ does not belong to $\supp L$. If $K$ is a connected 
component,  containing $h(z_L)$, of the complement of $\supp L$ in $\Bbb C$ 
and $h(z)\in K$ for some point $z\in\cl L$ then belongs to the midrib of $L$.
\endproclaim
\demo{Proof} We fix a neighbourhood $W$ of $z_L$ and a point $z\in\cl L$ such 
that $h(z)\in K$. Let $W^*=h(W)$. Since $h(z_L)\ne h(w)$ for all $w\ne z_L$ 
we can find a neighbourhood $V\Sbs W$ of $h(z_L)$ such that 
$h^{-1}(V)\Sbs W$. Let us take an approximation $g$ of $h$ by functions 
holomorphic on $\cl L$ such that $g(z)$ belongs to 
the connected component $K'$, containing $g(z_L)$, of the complement of 
the set $g(\supp L)$ in $\Bbb C$ and $g^{-1}(V)\Sbs W$. We consider on $K'$ 
the function $$u(\zeta)=\om(\zeta,g(W),K')$$ which is, evidently, greater 
than 0 at $g(z)$. 
Since $$\om(z,W,L)\ge u(g(z))>0$$ the point $z$ belongs to the midrib.\qed
\enddemo
It follows from this theorem that the leaf generated by the example of Wermer 
\cite{14} coincides with its midrib. 
\heading \S 5. Polynomial hulls\endheading
Let $D$ be a domain in $\Bbb C^n$ and let $\psi$ be an upper semicontinuous 
function on $D$. 
We consider the function 
$$P\psi=u(z)=\inf\frac1{2\pi}\il 0^{2\pi}\psi(f(e^{i\th}))\,d\th,$$where the 
infimum is taken over all mappings $f\in A(\ovr U,D)$ with $f(0)=z$. It was 
proved in \cite{8} that $u(z)$ is plurisubharmonic and  
$$u(z)=\sup\{v(z):v\le\psi \text{ and } v\text{ is psh 
}\}.$$\par 
Let $E\subset D$ be an open set. Then  the characteristic function $\chi_E$ 
of $E$ is lower semicontinuous and we define the {\it pluri-harmonic 
measure } $\om(z,E,D)$ of the set $E$ with respect to $D$ to be equal to 
$P(-\chi_E)$. If $F$ is an arbitrary set in $D$ then we define the 
pluri-harmonic measure $\om(z,F,D)$ to be equal to $\sup \om(z,E,D)$, 
where the supremum is taken over all open sets $E$ containing $F$. This 
function may be not upper  semicontinuous, so we take its regularization 
$$\om^*(z,F,D)=\varlimsup\limits_{w\to z}\om(z,F,D),$$ which  is 
plurisubharmonic.\par
If $X$ is a compact in $\Bbb C^n$ then we shall denote by $\hat X$ the 
polynomial hull of $X$. It follows from Theorem 3.3 that if a leaf $L$ is 
perfect and $X=\supp L$ then $\cl L\sbs\hat X$.\par
If $X\sbs D$ and $D$ is a 
Runge domain then $\hat K$ coincides with the holomorphic envelope of $K$ in 
$D$. The following theorem was proved in \cite{9}. For the reader's 
convenience we shall supply the proof.
\proclaim{Theorem 5.1} Let $D$ be a Runge domain and $K$  be a compact in 
$D$. Then $z_0\in\hat K$ if and only if for  every open set $E\supset K$ the 
pluri-harmonic measure $\om(z_0,E,D)=-1$.\endproclaim 
\demo{Proof} Let $z_0\in\hat K$ and $E$ is an open set containing $K$. Then 
there is a decreasing sequence of continuous psh functions $v_k(z)$ on $D$ 
converging to $\om(z,E,D)$ pointwise on $D$ \cite{11}. By \cite{2} for every 
$\eps>0$ there are holomorphic functions $f_{jk},j=1,2,\dots,N_k$ on $D$ and 
positive numbers $a_{jk}$ such that $$\max\limits_j\{a_{jk}\ln|f_{jk}|\}\le 
v_k\le\max\limits_j\{a_{jk}\ln|f_{jk}|\}+\eps$$ on $K$. Since $D$ is a Runge 
domain we may assume that all $f_{jk}$ are polynomials. We can take $k$ so 
big that $v_k\le-1+\eps$ on $K$. Then 
$v_k(z_0)\le-1+2\eps$ and this means that $\om(z_0,E,D)=-1$.\par 
Conversely, let $g$ be a polynomial. For every $\eps>0$ we consider the open 
set $E=\{z\in D:\,|g(z)|<\sup\limits_K|g(w)|+\eps\}\sps K$. If 
$\om(z_0,E,D)=-1$ then there is a  mapping $f\in A(\ovr U,D)$ with $f(0)=z_0$ 
such that $m\{\th:\,f(e^{i\th})\in E\}>2\pi(1-\eps)$. Therefore,
$$|g(z_0)|\le\frac1{2\pi}\int\limits_0^{2\pi}|g(f(e^{i\th}))|\, d\th\le\sup
\limits_K|g(w)|+M\eps, $$where $M=\sup\limits_D|g(w)|$. Hence, 
$z_0\in \hat K$.\qed\enddemo
The following theorem establishes something like analytic structure in 
$\hat X$.
\proclaim{Theorem 5.2} If $X$ is a compact in $\Bbb C^n$ and 
$z\in\hat X\sm X$ then there is a perfect leaf $L=\{f_j\}$ such that:
\roster\item 
$z=z_L$;\item $\supp L\sbs X$;\item $\cl L\sbs\hat X$;\item for every 
neighbourhood $V$ of $\hat X$ there is an integer $k$ such that analytic 
disks $f_j(U)$ are in $V$ when $j>k$.\endroster\endproclaim
\demo{Proof} By Theorem 5.1 we can find open neighbourhoods $W_j$ of $X$ 
with the intersection equal to $X$ and analytic disks $f_j$ such that 
the measure  of points $\zeta\in S$ with $f_j(\zeta)\in W_j$ is greater than 
$2\pi-j^{-1}$ and $f_j(0)=z$. Taking, if necessary, a subsequence we may 
assume that the sequence $L=\{f_j\}$ is a leaf. Evidently, $\supp L\sbs X$. 
We may also assume that $L$ is perfect, otherwise we replace $L$ by its 
perfection.\par
By the construction, $z=z_L$ and $\supp L\sbs X$. Since $\cl L$ belongs 
to the polynomial hull of $\supp L$, we see that $\cl L\sbs\hat X$. 
Therefore, for every neighbourhood $V$ of $\hat X$ there is an integer $k$ 
such that analytic disks $f_j(U)$ are in $V$ when $j>k$.\qed\enddemo
\heading\S 6. Maximal functions\endheading
The plurisubharmonic (psh) function $u$ on a domain $D\sbs\Bbb C^n$ is called 
{\it maximal } if for every domain $G\Sbs D$ and every psh function $v$ on 
$D$ the inequality $v\le u$ on $\bd G$ implies the inequality $v\le u$ on 
$G$. Bedford and Kalka proved in \cite{1} that if a maximal function $u$ is 
in $C^2(D)$ then through every point $z\in D$ we may draw a holomorphic disk 
such that the restriction of $u$ to this disk is harmonic. The set of all 
this disks is called {\it the Monge--Amp\'ere foliation for } $u$.\par
It was proved in \cite{9, Lemma 2.2.3, Theorem 2.2.1} that maximal functions 
on strongly pseudoconvex domains can be obtained as solutions of variational 
problems for mappings of the unit disk. More precisely:
\proclaim{Theorem 6.1} Let $D$ be a strongly pseudoconvex domain in 
$\Bbb C^n$ and let $\hvi<M$ be a continuous function on $\bd D$. We let 
$\hvi$ to be equal to $M$ on $D$. Let 
$$\Hvi(f)=\frac1{2\pi}\int\limits_0^{2\pi}\hvi(f(e^{i\th}))\,d\th$$ be a 
functional on the space $A(U,D)$ of all holomorphic mappings $f:U\to D$. 
Then the function $u(z)=P\hvi=\inf\Hvi(f)$, where the infimum is taken over 
all $f\in A(U,D)$ such that $f(0)=z$, is a maximal psh function on $D$, 
continuous up to the boundary, and equal $\hvi$ on $\bd D$.\par
Moreover, for every point $z\in D$ there is a sequence $\{f_j\}$ of 
holomorphic mappings $f_j:U\to D$ such that $f_j(0)=z$, $\supp\mu(dz,f_j)$ 
belongs to $\bd D$, and $\Hvi(f_j)\to u(z)$ as $j\to\infty$.\endproclaim
As the following example shows the existence of Monge--Amp\'ere foliations in 
the sense of Bedford and Kalka, is not valid in the general case.\par
\subheading{ Example 6.1 } Let $K$ be a compact in the boundary of the unit 
ball $B$, such that its polynomial hull $\hat K\sbs B$ doesn't contain 
analytic disks. (One can use examples of Stolzenberg and Wermer to build such 
$K$.) If $\hvi$ is a continuous negative function on $S=\bd B$ 
equal to $-1$ on $K$ and greater than $-1$ everywhere else, then we can find 
a maximal function $v$ which is continuous in $\ovr B$ and is equal to $\hvi$ 
on $S$. By Theorem 5.1 the function $v=P\hvi=-1$ only on $\hat K$. If 
$z\in\hat K$ and $f:U\to D$ is a holomorphic mapping such that $f(0)=z$ and 
$v\circ f$ is harmonic then $v\circ f\equiv-1$ and, hence, $f(U)\sbs\hat K$. 
But $\hat K$ does not contain analytic disks, so $f$ is trivial, which means  
that the Monge--Amp\'ere foliation for $v$ does not cover $\hat K$.\par
So to extend the result of Bedford and Kalka to the general case we have to 
look for more general constructions.\par
Let $L=\{f_j\}$ be a leaf such that $\supp L\sbs\bd D$ and let $G\sbs D$ be 
another domain, containing an  inner point $z$ of $L$. For a principal 
subleaf $L_z=\{g_j\}$  we consider a connected component $G_j$, containing 
the origin, of the set $g^{-1}_j(G)$ and let 
$h_j$ be a holomorphic covering mapping of $G_j$ by $U$ such that $h_j(0)=0$. 
Then {\it the restriction } of $L$ on $G$ at $z$ is the perfection of the 
sequence $\{g_j\circ h_j\}$. We say that an upper semicontinuous function $u$ 
on $\cl L$ is {\it maximal } if for all domains $G$ and $H$ such that 
$G\Sbs H\Sbs D$ and every plurisubharmonic (resp. plurisuperharmonic) 
function $v$ on $H$, the inequality $u\ge v$ (resp. $u\le v$) on the 
intersection of $\cl L$ and $\bd G$ implies 
that $u(z)\ge v(z)$ (resp. $u(z)\le v(z)$) on the intersection of $\cl L$ 
and $G$. 
\proclaim{Theorem 6.2} Let $D$ be a bounded domain in $\Bbb C^n$ and let $u$ 
be a continuous plurisubharmonic function on $D$. We consider the functional 
$$\Hvi(f)=\frac1{2\pi}\int\limits_0^{2\pi}u(f(e^{i\th}))\,d\th$$ on the set 
$A(\ovr U,D)$. Suppose that for a point $z_0\in D$ there is a leaf 
$M=\{f_j\}$ such that $f_j\in A(U,D)$, $z_0=z_M$, $\supp M\sbs\bd D$ and 
$$u(z_0)=\lim\limits_{j\to\infty}\Hvi(f_j).$$
Then $u$ is maximal on the perfection $L$ of $M$.\endproclaim
\demo{Proof} We may assume that $L=M$. It is clear that $u(z)\le \Hvi(f)$ for 
every $z\in D$ and for every $f\in A(U,D)$ such that $f(0)=z$. Let us take 
arbitrary domains $G\Sbs H\Sbs D$ and let $v$ be a plurisubharmonic function 
on $H$ such that $u\ge v$ on the intersection $K$ of $\cl L$ and $\bd G$. We 
may assume that both $u$ and $v$ are less than some number $A$ on $H$.\par
Suppose that there is a point $z$ in the intersection $F$ of $\cl L$ and $G$ 
such that $u(z)<v(z)$. We fix $\eps>0$ such that $u(z)<v(z)-3\eps$. Since $u$ 
is continuous we can find a neighbourhood $V$ of $K$ such that $u>v-\eps$ on 
$V$ and since $L$ is perfect we may assume that if $f_j(\zeta)\in\bd G$ then 
$f_j(\zeta)\in V$. For each $j$ we consider domains $G_j=f_j^{-1}(G)$ which 
are unions of connected disjoint domains $G_{jk}$.\par
Taking averages if necessary we may assume that $v$ is continuous so 
$u(z)<v(z)-\eps$ on a ball $B=B(z,r)$. We denote by $B_j$ sets $f_j^{-1}(B)$. 
The boundary $\Gm_{jk}$ of domains $G_{jk}$ consists of points of $U$ and 
points of $S$. We denote the first part by $\gm_{jk}$, and the second by 
$\gm'_{jk}$. If $u_j(\zeta)=u(f_j(\zeta))$ and $v_j(\zeta)=v(f_j(\zeta))$ 
then $u_j>v_j-\eps$ on all $\gm_{jk}$. Since $\supp L\sbs\bd D$ harmonic 
measures $a_j=\om(0,\gm'_{jk},U)$ converge to 0. If $E_j$ is the set of 
points in $B_j$ where $\om(\zeta,\gm'_{jk},G_{jk})>s>0$ then, evidently, 
$$\om(0,E_j,U)<\frac{a_j}s.$$ Let us introduce harmonic functions $h_j$ on 
$U$ with boundary values of $u_j$ on $S$. Since $h_j\ge u_j$ on $U$ the 
function $$h_j\ge v_j-\frac{Aa_j}s-\eps$$ on the set $C_j=B_j\sm E_j$ and, 
therefore, $$h_j>u_j-\frac{Aa_j}s+2\eps>u_j+\eps$$ on $C_j$ when $j$ is 
sufficiently big. Since the leaf $L$ is perfect, harmonic measures 
$$\om(0,B_j,U)>a>0$$ and, therefore, $$\om(0,C_j,U)>a-\frac{a_j}s>\frac a2$$ 
when $j$ is sufficiently big. So if we consider the universal holomorphic 
covering mappings $g_j$ of $U\sm C_j$ by $U$, and if we take mappings 
$q_j=f_j\circ g_j$ then $$u(z_0)\le\Hvi(q_j)<\Hvi(f_j)-\frac{a\eps}2,$$and we 
get the contradiction.\par
Now if $v$ is plurisuperharmonic and $v\ge u$ on $K$ then, as before, we may 
assume that $v$ is continuous on $H$. Repeating estimates from the first part 
of the theorem we may prove that for every $r>0$ and for every $\eps>0$ the 
harmonic measure of points in $B_j$ where $v_j\ge u_j$ is uniformly greater 
than 0 in $U$. For us it is enough to know that this set is non-empty because 
the combination of continuity of $u$ and $v$ and arbitrariness of $r$ and 
$\eps$ implies that $v(z)\ge u(z)$.
\qed\enddemo
The combination of Theorem 6.1 and 6.2 gives us 
\proclaim{Theorem 6.3} A continuous plurisubharmonic function $u$ on a 
strongly pseudoconvex domain 
$D$ in $\Bbb C^n$ is maximal if and only if through every point $z\in D$ 
passes the cluster of a perfect leaf $L$ such that:\roster\item $u$ is 
maximal on $\cl L$;\item $\supp L\sbs\bd L$.\endroster\endproclaim

\enddocument